\documentclass[10pt]{amsart}
\usepackage{}
\usepackage{color}
\usepackage{amssymb}

\newcommand{\ov}{\overline}

\newtheorem{theorem}{Theorem}[section]
\newtheorem{lemma}[theorem]{Lemma}

 \theoremstyle{definition}

\theoremstyle{remark}
\newtheorem{remark}[theorem]{Remark}

\numberwithin{equation}{section}

\begin{document}

\title[Second order estimates]
{second order estimates for complex Hessian equations on Hermitian manifolds}

\author{Weisong Dong}
\address{School of Mathematics, Tianjin University,
         Tianjin, P.R.China, 300354}
\email{dr.dong@tju.edu.cn}

\author{Chang Li}
\address{School of Mathematical Sciences, Peking University, Beijing, P.R.China, 100871}
\email{chang\_li@pku.edu.cn}

%\date{}

\begin{abstract}
%We consider the complex Hessian equations with both sides depending on gradients. We derive the second order estimates for $\chi$-plurisubharmonic solutions on compact Hermitian manifolds.
We derive second order estimates for $\chi$-plurisubharmonic solutions of complex Hessian equations with right hand side depending on the gradient on compact Hermitian manifolds.



\end{abstract}

\maketitle

\section{Introduction}

Let $(M, \omega)$ be a compact complex manifold of complex dimension $n \geq 2$.
For any smooth function $u \in C^{\infty}(M)$, let $\chi(z,u)$ be a smooth real (1,1) form on $M$
and $\psi(z, v, u) \in C^{\infty}\left(\left(T^{1,0}(M)\right)^{*} \times \mathbb{R}\right)$ be a positive function.
The following equation which we shall call it the complex Hessian equation:
\begin{equation}\label{Sk equation}
(\chi(z, u)+\sqrt{-1} \partial \overline{\partial} u)^{k} \wedge \omega^{n-k}=\psi(z, Du, u) \omega^{n},
\end{equation}
for $1\leq k \leq n$,
can be viewed as an intermediate equation between the Laplace equation ($k=1$) and the complex Monge-Amp\`ere equation ($k=n$),
where $D$ is the covariant derivative with respect to the given metric $\omega$.

A function $u \in C^2(M)$ is called admissible if
$g=\chi (z, u)+\sqrt{-1} \partial \overline{\partial} u \in \Gamma_{k}(M)$.
We formally define $\Gamma_k (M)$ in \eqref{Gamma-k}.
In particular, $u$ is called $\chi$-plurisubharmonic  if $g=\chi (z, u)+\sqrt{-1} \partial \overline{\partial} u > 0$ (i.e. $g\in \Gamma_n (M)$).
When $k=1$, \eqref{Sk equation} is just a quasilinear equation which is well understood.
Otherwise, it is fully nonlinear and a natural approach to solve \eqref{Sk equation} is the continuity method which reduces the solvability to a priori estimates of solutions up to the second order. Higher order
estimates follow from Evans-Krylov theory and Schauder estimate, which was done by Tosatti, Wang, Weinkove and Yang \cite{T-W-W-Y}.
Note that Tian \cite{Tian} presented a new proof of the $C^{2,\alpha}$ estimates (for real and complex Monge-Amp\`ere equations), which not only weakens the regularity assumptions on $\psi(z)$ but also can be applied to more general nonlinear elliptic systems. In \cite{Tian2}, Tian extended his method to the conic case.
In this paper, we are mainly concerned with the second order estimate of \eqref{Sk equation}.

On compact K\"ahler manifolds, the existence and regularity theory for complex Monge-Amp\`ere equations has been studied for a long time since Yau proved the Calabi conjecture in \cite{Yau},
where he also studied the case when the right-hand side may degenerate or have poles.
Later after that, Cheng-Yau \cite{C-Y}, Kobayashi \cite{R.K} and Tian-Yau \cite{T-Y} \cite{T-Y2} generalized the conjecture to complete noncompact K\"ahler manifolds.
We refer the readers to a survey on complex Monge-Amp\`ere equations by Phong-Song-Sturm \cite{P-S-S}.
The complex Hessian equation
was first considered on domains in $\mathbb{C}^n$ by Li \cite{S.Y.Li}, then by Blocki \cite{Blocki}.
On compact K\"ahler manifolds, Hou \cite{Z.L.H} and Jbilou \cite{A.J} proved the existence of a smooth admissible solution of \eqref{Sk equation}
independently when $\chi= \omega$ and
$\psi = \psi (z)$ by assuming the nonnegativity of the holomorphic bisectional curvature of $\omega$.
This curvature assumption was removed in the work of Hou-Ma-Wu \cite{H-M-W} and Dinew-Ko\l odziej \cite{D-K}.
The complex Hessian equation also appears in many geometric problems, such as the $J$-flow studied by Song-Weinkove \cite{S-W} and quaternionic geometry by Alekser-Verbitsky \cite{A-V}.

There has been growing interest in extending the above results to non-K\"ahler settings.
For the complex Monge-Amp\`ere equation, it was solved by the work of Cherrier \cite{Cherrier} and Tosatti-Weinkove \cite{T-W1} \cite{T-W2} on Hermitian manifolds,
as well as by Chu-Tosatti-Weinkove \cite{C-T-W} on almost Hermitian manifolds.
The related Dirichlet problem was studied by Guan-Li \cite{G-L}.
Zhang \cite{D.K.Zhang} solved the equation \eqref{Sk equation} on compact Hermitian manifolds when $\chi=\omega$ and $\psi=\psi(z)$.
Chu, Huang, and Zhu \cite{C-H-Z3}
obtained the second order estimate with $k=2$ on an almost Hermitian manifold.
See also Sun \cite{Sun1} \cite{Sun2} and Sz\'ekelyhidi \cite{G.S} for some other equations.

When $\psi=\psi(z, Du, u)$, \eqref{Sk equation} has been much less studied so far.
Actually, it is still an open problem even for the real counterparts of \eqref{Sk equation}
to derive a second order estimate.
Recently, Guan-Ren-Wang \cite{G-R-W} solved this problem for convex solutions.
Borrowing the idea in \cite{G-R-W} and adapting the techniques developed by Li-Ren-Wang \cite{L-R-W} for real Hessian equations,
Phong, Picard and Zhang \cite{P-P-Z C^2} obtained the second order estimate for $\chi$-plurisubharmonic solutions of \eqref{Sk equation} on K\"ahler manifolds.
In this paper,
motivated by \cite{P-P-Z C^2},
we study \eqref{Sk equation} on compact Hermitian manifolds
and derive the second order estimate for $\chi$-plurisubharmonic solutions.
Our main result is stated as follows.
\begin{theorem}\label{thm 1}
Let $(M,\omega)$ be a compact Hermitian manifold of complex dimension n.
Suppose $u \in C^{4}(M)$ is a $\chi$-plurisubharmonic solution
of \eqref{Sk equation} and $\chi(z, u) \geq \varepsilon \omega$.
Then we have the uniform second order derivative estimate
\begin{equation}
|D \ov{D} u|_{\omega} \leq C,
\end{equation}
where $C$ is a uniform constant depending only on $(M,\omega)$, $\varepsilon$, $n$, $k$, $\chi$, $\psi$, $\sup_{M}|u|$, $\sup_{M}|D u|$.
\end{theorem}

\begin{remark}
The above estimate can be stated for $g\in \Gamma_{k+1}$, see the remark after Theorem 1 and Remark 2 in \cite{P-P-Z C^2}.
An interesting question is to derive the estimate with the more natural elliptic condition $g\in \Gamma_k$.
If $k=n$, $\chi$-plurisubharmonic is already the natural assumption for ellipticity of \eqref{Sk equation}.
So our result generalizes the estimate for Monge-Amp\`ere equations on compact Hermitian manifolds.
If $k=2$, the estimate was derived without the plurisubharmonicity assumption, see \cite{C-H-Z3}.
\end{remark}

An important case of \eqref{Sk equation} with $\psi$ depending on $Du$ is the Fu-Yau equation.
As a reduced version of generalized Strominger system in higher dimensions, Fu and Yau in \cite{F-Y2} introduced a fully nonlinear equation which can be rewritten as a $\sigma _2$-type equation with specific right hand side $\psi (z, Du, u)$.
When $n=2$, the Fu-Yau equation is equivalent to the Strominger system on a toric fibration over a $K3$ surface constructed by Goldstein and Prokushki \cite{G-P}, which was solved by Fu and Yau in \cite{F-Y2, F-Y1}.
For higher dimension $n$, the corresponding problem on compact K\"ahler or Hermitian manifolds has been well studied,
see \cite{P-P-Z1, P-P-Z2, C-H-Z1, P-P-Z3, C-H-Z2}.

The dependence on the gradient of $u$ in \eqref{Sk equation} creates substantial new difficulties, due to the appearance of terms such as $|D D u|^{2}$ and $|D \ov{D} u|^{2}$
when one differentiates the equation twice.
A consequence of this is that we cannot control the bad third order terms straightly as in Guan-Jiao \cite{G-J}.
Furthermore, it is made more difficult by the differences between the real case and the complex case to control the negative third order terms due to complex conjugacy.
The differences can be seen right through the definition of third order terms in Section 3 that the good terms $``B"$ and $``D"$
compared with those in \cite{G-R-W}
get doubled in the real case.
These are the main difficulties that have been overcome in \cite{P-P-Z C^2} on K\"ahler manifolds.
In the Hermitian case, there are more bad terms from the appearance of the form $T*D^3u$, where $T$ is the torsion of $\omega$ and $D^3u$ represents the third derivatives of the solution $u$. To control these terms, we modify the auxiliary function. This gives us a little more good third order terms which are sufficient to push the argument through.

Many authors are attracted by some other fully nonlinear equations involving gradient terms.
In contrast to the equation \eqref{Sk equation}, those gradient terms appear together with the Hessian of the solution,
that is they are in $\chi$.
For exmaple, some related equations were solved by Sz\'ekelyhidi-Tosatti-Weinkove in their work on the Gauduchon conjecture \cite{S-T-W},
and see Guan-Nie \cite{GN} for related work.
Also, see a recent work by Tosatti-Weinkove \cite{T-W19} where they solved the complex Monge-Amp\`ere-type equation with gradient terms
left open by Yuan \cite{Yuan}.
Very recently, Yuan informed the authors that he solved the same problem as in \cite{T-W19}.
Moreover, he studied the Dirichlet problem of more general fully nonlinear equations on Hermitian manifolds in \cite{Yuan2}.
When $\chi$ may depend on the gradient of $u$ in \eqref{Sk equation},
Zhang \cite{X.W.Zhang} studied the gradient estimate and
Feng-Ge-Zheng \cite{Feng-Ge-Zheng} derived the second order estimate.
Our method here can also deal with such extra gradient terms. Suppose $a$ is a smooth $(1,0)$-form on $M$.
Let
\begin{equation}
\tilde{g}=\tilde{\chi}(z, u, Du) +\sqrt{-1} \partial \overline{\partial} u,
\end{equation}
where $\tilde{\chi}(z, u, Du)=\chi(z, u)+\sqrt{-1} a \wedge \overline{\partial} u-\sqrt{-1} \overline{a} \wedge \partial u$,
and consider the equation
\begin{equation}\label{Sk equation 2}
\tilde{g}^{k} \wedge \omega^{n-k}=\psi(z, Du, u) \omega^{n},
\end{equation}
for $1\leq k \leq n$. In analogy to Theorem \ref{thm 1}, we have the following result.

\begin{theorem}\label{thm 2}
Let $(M,\omega)$ be a compact Hermitian manifold of complex dimension n.
Suppose $u \in C^{4}(M)$ is a $\tilde{\chi}$-plurisubharmonic solution
of \eqref{Sk equation 2} and $\chi(z, u) \geq \varepsilon \omega$.
Then we have the uniform second order derivative estimate
\begin{equation}
|D \ov{D} u|_{\omega} \leq C,
\end{equation}
where $C$ is a uniform constant depending only on $(M,\omega)$, $\varepsilon$, $n$, $k$, $\chi$, $\psi$, $a$, $\sup_{M}|u|$, $\sup_{M}|D u|$.
\end{theorem}

\begin{remark}
If $\psi$ does not depend on $Du$, one can get the above estimate
for admissible solutions
by estimating the largest eigenvalue of $\tilde g$ as in \cite{G.S}.
\end{remark}

The rest of the paper is organized as follows.
In Section 2, we introduce some notations and calculations.
We prove Theorem \ref{thm 1} in Section 3  and give an outline of the proof of Theorem \ref{thm 2} in Section 4.

\textbf{Acknowledgements}:
The second named author would like to thank his advisor Prof. Gang Tian for constant encouragement and support.
This paper was done while both authors were visiting Rutgers, The State University of New Jersey.
Both authors would like to thank the department of mathematics for the warm hospitality.
The first named author's visiting is supported by the China Scholarship Council (File No.201806255014).
The first named author is also partially supported by a NSFC grant No.11801405.

\section{Preliminaries}

To be more clear, we shall rewrite the equation \eqref{Sk equation} in local coordinates.
First, we introduce some notations.
For  $\lambda = \left(\lambda_{1}, \cdots, \lambda_{n}\right) \in \mathbb{R}^{n}$,
the $k$-th elementary symmetric function is defined by
$\sigma_{k}(\lambda)=\sum \lambda_{i_{1}} \lambda_{i_{2}} \cdots \lambda_{i_{k}}$,
where the sum is over $\{1\leq i_{1}<\cdots<i_{k}\leq n\}$.
Define the cone $\Gamma_{k}$ by
\begin{equation}
\Gamma_{k} :=\left\{\lambda \in \mathbb{R}^{n} : \sigma_{j}(\lambda)>0, j=1, \cdots, k\right\}.
\end{equation}
It was shown in \cite{C-N-S} that $\Gamma_k$ is an open convex symmetric cone in $\mathbb{R}^n$.
Let $\lambda (a_{i \overline{j}})$ denote the eigenvalues of a Hermitian symmetric matrix $(a_{i \overline{j}})$. Define
$\sigma_{k}(a_{i\overline{j}})=\sigma_{k}(\lambda(a_{i \overline{j}}))$.
This definition can be naturally extended to complex manifolds.
Let $A^{1,1}(M)$ be the space of smooth real $(1,1)$-forms on $(M, \omega)$. For any $h \in A^{1,1}(M)$, define
\begin{equation}
\sigma_{k}(h)=\left( \begin{array}{l}{n} \\ {k}\end{array}\right) \frac{h^{k} \wedge \omega^{n-k}}{\omega^{n}}.
\end{equation}
Now, $\Gamma_{k}$ can be defined on $M$ as follows
\begin{equation}
\label{Gamma-k}
\Gamma_{k}(M) :=\left\{h \in A^{1,1}\left(M, \mathbb{R}^{n}\right) : \sigma_{j}(h)>0, j=1, \cdots, k\right\}.
\end{equation}
For a function $u \in C^{\infty}(M)$, we denote
\begin{equation}
g=\chi(z, u)+\sqrt{-1} \partial \overline{\partial} u.
\end{equation}
With the above notation, in local coordinates, \eqref{Sk equation} can be rewritten as follows:
\begin{equation}\label{Rewrite Sk equ}
\sigma_{k}(g)=\sigma_{k}\left(\chi_{i \overline{j}}+u_{i \overline{j}}\right)=\psi(z, D u, u).
\end{equation}

Recall that if $F(A)= f(\lambda_1, \ldots, \lambda_n)$ is a symmetric function of the eigenvalues of a Hermitian matrix $A=(a_{i \ov j})$, then at a diagonal matrix $A$ with distinct eigenvalues, we have (see \cite{J.M.B}),
\begin{align}
\label{symmetric func 1th deriv} F^{i\ov j} =&\ \delta_{ij} f_i,\\
\label{symmetric func 2th deriv} F^{i\ov j, r\ov s} w_{i\ov j k} w_{r\ov s \ov k} =&\ \sum f_{ij} w_{i\ov i k} w_{j \ov j \ov k} + \sum_{p\neq q}\frac{f_p - f_q}{\lambda_p-\lambda_q} | w_{p\ov q k}|^2,
\end{align}
where $F^{i\ov j}=\frac{\partial F}{\partial a_{i \ov j}}$,
$F^{i\ov j,r\ov s}=\frac{\partial^2 F}{\partial a_{i \ov j}\partial a_{r \ov s}}$, and $w_{i\ov jk}$ is an arbitrary tensor.

In local complex coordinates $\left(z_{1}, \ldots, z_{n}\right)$, the subscripts of a function $u$ always denote the covariant derivatives of $u$ with respect to $\omega$ in the directions of the local frame $\left(\partial / \partial z^{1}, \ldots, \partial / \partial z^{n}\right)$. Namely,
\begin{equation}\nonumber
u_{i}=D_{i}u=D_{\partial / \partial z^{i}} u, \; u_{i \overline{j}}=D_{\partial / \partial \overline{z}^{j}} D_{\partial / \partial z^{i}} u, \; u_{i \overline{j} l}=D_{\partial / \partial z^{l}} D_{\partial / \partial \overline{z}^{j}} D_{\partial / \partial z^{i}} u.
\end{equation}
We have the following commutation formula on Hermitian manifolds (see \cite{T-W3} for more details):
\begin{equation}\label{order}
\begin{aligned}
u_{i \overline{j} \ell}= &\  u_{i \ell \overline{j}}-u_{p} R_{\ell \overline{j} i}{}^p, \;\;\\
u_{p \overline{j} \overline{m}}= &\  u_{p \overline{m} \overline{j}}-\overline{T_{m j}^{q}} u_{p \overline{q}},\;\;\\
u_{i \overline{q} \ell}= &\  u_{\ell \ov q i}-T_{\ell i}^{p} u_{p \overline{q}},
\end{aligned}
\end{equation}
\begin{equation}
u_{i \overline{j} \ell \overline{m}} =  u_{\ell \overline{m} i \overline{j}}
+u_{p \overline{j}} R_{\ell \overline{m} i}{}^{p}
 - u_{p \overline{m}} R_{i \overline{j} \ell}{}^{p}
- T_{\ell i}^{p} u_{p \overline{m} \overline{j}}
-\overline{T_{m j}^{q}} u_{\ell \overline{q} i}
 -T_{i \ell}^{p} \overline{T_{m j}^{q}} u_{p \overline{q}}.
\end{equation}

We use the notation
\begin{equation}
\sigma_{k}^{p \overline{q}}=\frac{\partial \sigma_{k}(g)}{\partial g_{p \overline{q}}}, \quad \sigma_{k}^{p \overline{q}, r \overline{s}}=\frac{\partial^{2} \sigma_{k}(g)}{\partial g_{p \overline{q}} \partial g_{r \overline{s}}}, \quad
\mathcal{F}=\sum_{p} \sigma_{k}^{p \overline{p}}.
\end{equation}
We also use the following notation as in \cite{P-P-Z C^2}
\begin{equation}
\begin{aligned}
|D D u|_{\sigma \omega}^{2}=\sigma_{k}^{p \overline{q}} \omega^{m \overline{\ell}} u_{mp} u_{\ov \ell \ov q},\;
|D \ov{D} u|_{\sigma \omega}^{2}=\sigma_{k}^{p \overline{q}} \omega^{m \overline{\ell}} u_{p\ov\ell} u_{m\ov q}.
\end{aligned}
\end{equation}
and
\begin{equation}
|\eta|_{\sigma}^{2}=\sigma_{k}^{p \overline{q}} \eta_{p} \eta_{\overline{q}},
\end{equation}
for any 1-form $\eta$.

Now we do some basic calculations which are used in next Section.
In the following, $C$ will be a uniform constant depending on the known data as in Theorem \ref{thm 1}, but may change from line to line.

Our calculations are carried out at a point $z$ on the manifold $M$, and we use coordinates such that at this point $\omega=\sqrt{-1} \sum \delta_{k \ell} dz^{k} \wedge d \overline{z}^{\ell}$ and $g_{i \ov{j}}$ is diagonal.
Differentiating \eqref{Rewrite Sk equ} yields
\begin{equation}\label{differential equ}
\sigma_{k}^{p \overline{q}} D_{i} g_{p \overline{q}}=D_{i} \psi.
\end{equation}
Differentiating the equation a second time gives
\begin{equation}\label{diff equ second time}
\begin{split}
&\sigma_{k}^{p \overline{q}} D_{\ov{j}} D_{i} g_{p \ov q}+\sigma_{k}^{p \overline{q}, r \overline{s}}D_{\overline{j}} g_{r \overline{s}} D_{i} g_{p \overline{q}}=D_{\ov j}D_{i} \psi \\
\geq & -C(1+|DDu|^{2}+|D \ov{D} u|^{2})+\sum_\ell \psi_{v_\ell} u_{\ell i \ov{j}} + \sum_\ell \psi_{\ov{v}_\ell} u_{\ov{\ell} i \ov{j}}\\
\geq & -C(1+|DDu|^{2}+|D \ov{D} u|^{2})+\sum_\ell \psi_{v_\ell} g_{i \ov{j} \ell} + \sum_\ell \psi_{\ov{v}_\ell} g_{i \ov{j} \ov{\ell}}-C \lambda_{1}.
\end{split}
\end{equation}
Direct calculation gives the estimate
\begin{align}\label{4th order term}
\sigma_k^{p \ov{q}} D_{\ov{q}} D_p g_{i \ov{j}}
\geq \sigma_k^{p \ov{q}} D_{\ov{q}} D_p D_{\ov{j}} D_i u - C(1 + \lambda_1) \mathcal{F}.
\end{align}
Commuting derivatives yields that
\begin{equation}
\begin{aligned}\label{commute 4th order deri}
D_{\ov{q}} D_p D_{\ov{j}} D_i u =&\  D_{\ov{j}} D_i D_{\ov{q}} D_p u - R_{i \ov{j} p}{}^{a} u_{a \ov{q}} + R_{p \ov{q} i}{}^{a} u_{a \ov{j}} \\
&\ -T_{p {i}}^{a} u_{a \ov{q} \ov{j}}-\ov{T_{q j}^{a}} u_{p \ov a i} -T_{i p}^{a} \ov{T_{q j}^{b}} u_{a \ov{b}}\\
=&\ D_{\ov{j}} D_i g_{p \ov{q}} - D_{\ov{j}} D_i \chi_{p \ov{q}} -  R_{i \ov{j} p}{}^{a} u_{a \ov{q}} + R_{p \ov{q} i}{}^{a} u_{a \ov{j}} \\
&\ -T_{p {i}}^{a} u_{a \ov{q} \ov{j}} -\ov{T_{q j}^{a}} u_{p \ov a i}-T_{i p}^{a} \ov{T_{q j}^{b}} u_{a \ov{b}}.
\end{aligned}
\end{equation}
Therefore, by \eqref{diff equ second time}, \eqref{4th order term}, \eqref{commute 4th order deri}, we see
\begin{equation}
\begin{aligned}\label{4th order term after using equ}
& \sigma_{k}^{p \overline{q}} D_{\overline{q}} D_{p} g_{i \overline{j}}\\
 \geq &\ -\sigma_{k}^{p \overline{q}, r \overline{s}} D_{\ov{j}} g_{r \overline{s}} D_{i} g_{p \overline{q}}+\sum_\ell \psi_{v_\ell} g_{i \ov{j} \ell} + \sum_\ell \psi_{\ov{v}_\ell} g_{i \ov{j} \ov{\ell}}-\sigma_{k}^{p \ov{q}} (T_{p i}^{a} u_{a \overline{q} \overline{j}} \\
&\ +\overline{T_{q j}^{a}} u_{p \overline{a} i})-C(1+|DDu|^{2}+|D \ov{D}u|^{2}+\mathcal{F}+\lambda_{1} \mathcal{F}+\lambda_{1}).
\end{aligned}
\end{equation}
Then,
\begin{equation}
\begin{aligned}\label{(Du)_pq}
& \sigma_{k}^{p \overline{q}}|D u|_{p \overline{q}}^{2}\\
=&\ \sigma_{k}^{p \overline{q}}\left(u_{mp\overline{q}}D^m u+u_{m} u_{\overline{\ell} p \overline{q} } \omega ^{m \overline{\ell}}\right)+|DDu|_{\sigma \omega}^{2}+|D \ov{D} u|_{\sigma \omega}^{2} \\
=&\ \sigma_{k}^{p \overline{q}}\left(u_{p \ov q m}+(T_{mp}^{t}u_{t})_{\ov q}+u_{t}R_{m \ov q p}{}^{t}\right) D^m u \\
&\ +\sigma_{k}^{p \overline{q}} u_{m} \omega^{m \ov \ell}(u_{p \ov q \ov \ell}-\ov{T_{q \ell}^{t}}u_{p \ov t})+|DDu|_{\sigma \omega}^{2}+|D \ov{D} u|_{\sigma \omega}^{2} \\
=&\ \sigma_{k}^{p \overline{q}} D_{m}\left(g_{p \overline{q}}-\chi_{p \overline{q}}\right) D^{m} u+\sigma_{k}^{p \overline{q}} (T_{mp}^{t})_{\ov q} u_{t} D^{m}u+\sigma_{k}^{p \overline{q}} T_{mp}^{t} u_{t \ov{q}} D^{m}u \\
& \ +\sigma_{k}^{p \overline{q}} u_{t} R_{m \ov q p}{}^{t} D^m u+\sigma_{k}^{p \overline{q}} u_m \omega^{m \ov \ell} D_{\ov \ell}\left(g_{p \overline{q}}-\chi_{p \overline{q}}\right) \\
&\ - \sigma_{k}^{p \overline{q}} u_m \omega^{m \ov{\ell}} \ov {T_{q \ell}^{t}}u_{p \ov t}
+|D D u|_{\sigma \omega}^{2}+|D \ov{D} u|_{\sigma \omega}^{2}.
\end{aligned}
\end{equation}
Using the Cauchy inequality, we have
\begin{align}
|\sigma_{k}^{p \overline{q}} T_{mp}^{t} u_{t \ov{q}} D^{m}u|
 \leq \sum_{p,t} \sigma_{k}^{p \overline{p}}\left(\frac{1}{4}|u_{t \ov p}|^2+(C_p ^t)^2\right) \nonumber
\leq \frac{1}{4}|D \ov{D} u|_{\sigma \omega}^{2}+C\mathcal{F},
\end{align}
where $C_p ^t=T_{mp}^{t} D^{m}u$.
Similarly, we have
\begin{align}
|-\sigma_{k}^{p \overline{q}} u_{m} \omega^{m \overline{\ell}} \overline{T_{q \ell}^{t}} u_{p \overline{t}}| \leq \frac{1}{4}|D \ov{D} u|_{\sigma \omega}^{2}+C\mathcal{F} \nonumber.
\end{align}
Substituting the above two inequalities in \eqref{(Du)_pq}, we get
\begin{equation}
\begin{aligned}
&\ \sigma_{k}^{p \overline{q}}|D u|_{p \overline{q}}^{2} \\
\geq &\ \sigma_{k}^{p \overline{q}} D_{m} g_{p \overline{q}} D^{m} u+\sigma_{k}^{p \overline{q}} u_m  D^m g_{p \overline{q}}+|D D u|_{\sigma \omega}^{2}+\frac{1}{2}|D \ov{D} u|_{\sigma \omega}^{2}-C\mathcal{F} \\
= &\  D_{m}(\sigma_k) u_{\ov{\ell}}\omega^{m \ov{\ell}}+D_{\ov{\ell}}(\sigma_k)u_m \omega^{m \ov{\ell}}+|DDu|_{\sigma \omega}^2+\frac{1}{2}|D\ov{D} u|_{\sigma \omega}^2-C\mathcal{F}.
\end{aligned}
\end{equation}
Using the differential equation \eqref{differential equ}, we obtain
\begin{equation}
\begin{aligned}\label{derivative of Du}
\sigma_{k}^{p \overline{q}}|D u|_{p \overline{q}}^{2}
\geq &\; 2 \operatorname{Re}\left\{\sum_{p, m}\left(D_{p} D_{m} u D_{\overline{p}} u+D_{p} u D_{\overline{p}} D_{m} u\right) \psi_{v_{m}}\right\}\\
&\; -C-C\mathcal{F}+|D D u|_{\sigma \omega}^{2}+\frac{1}{2}|D \ov{D} u|_{\sigma \omega}^{2}.
\end{aligned}
\end{equation}
We also compute that
\begin{equation}\label{derivative of u}
-\sigma_{k}^{p \overline{q}} u_{p \overline{q}}=\sigma_{k}^{p \overline{q}}\left(\chi_{p \overline{q}}-g_{p \overline{q}}\right) \geq \varepsilon \mathcal{F}-k \psi.
\end{equation}

\section{Proof of Theorem \ref{thm 1}}

We denote by $\lambda_{1}, \lambda_{2}, \dots, \lambda_{n}$ the eigenvalues of $g_{i \overline{j}} =\chi_{i \ov{j}}+u_{i \overline{j}}$ with respect to $\omega$.
When $k=1$, the equation \eqref{Sk equation} becomes
\begin{equation}
\Delta_{\omega} u+\operatorname{Tr}_{\omega} \chi(z, u)=n \psi(z, D u, u).
\end{equation}
It follows that $\Delta_\omega u$ is bounded, and the desired estimate follows in turn from the positivity of $g$. Henceforth, we assume that $k \geq 2$.

We apply the maximum principle to the following test function:
\begin{equation}\label{test}
G=\log P_{m}+\varphi(|Du|^{2})+\phi(u),
\end{equation}
where $P_{m}=\sum_{j} \lambda_{j}^{m}$. Here, $\varphi$ and $\phi$ are positive functions to be determined later,
which satisfy the following assumptions
\begin{equation}
\label{testfun}
\varphi''-\phi''(\frac{\varphi'}{\phi'})^{2} \geq 0, \;\varphi' > 0, \; \phi' < 0.
\end{equation}
We assume that the maximum of $G$ is achieved at some point $p \in M$.
We choose the coordinate system centered at $p$ such that $\omega=\sqrt{-1} \sum \delta_{k \ell} dz^{k} \wedge d \overline{z}^{\ell}$ and $g_{i \ov{j}}$ is diagonal
with the ordering $\lambda_{1} \geq \lambda_{2} \geq \cdots \geq \lambda_{n}>0$.

Differentiating $G$ ,
we first obtain the critical equation
\begin{equation}\label{critical euq}
\frac{D P_{m}}{P_{m}}+\varphi' D|D u|^{2}+\phi' Du=0.
\end{equation}
Differentiating $G$ a second time, using \eqref{symmetric func 2th deriv} and contracting with $\sigma_{k}^{p \overline{q}}$ yields
\begin{equation}\label{diffrrentiate G second time}
\begin{aligned}
0 \geq &\ \frac{m}{P_m}  \sum_j \lambda_j^{m-1} \sigma_k^{p \ov{p}} D_{\ov p} D_{p} g_{j \ov{j}}+\frac{m\sigma_k^{p \ov{p}}}{P_m} (m-1)\sum_j \lambda_j^{m-2} |D_p g_{j \ov{j}}|^2\\
&\ +\frac{m\sigma_k^{p \ov{p}}}{P_m} \sum_{i \neq j} \frac{\lambda_i^{m-1} - \lambda_j^{m-1}}{\lambda_i - \lambda_j} |D_p g_{i \ov{j}}|^2 +\sigma_k^{p \ov{p}} \left( \phi'' D_p u D_{\ov p}u+\phi' u_{p \ov p} \right)\\
&\ + \sigma_k^{p \ov{p}}\left( \varphi'' D_p |Du|^2 D_{\ov p}|Du|^2+\varphi'|Du|_{p \ov p}^2 \right)- \frac{|D P_m|^2_\sigma}{P_m^2}.
\end{aligned}
\end{equation}
Here we used the notation introduced in Section 2.

Using the critical equation \eqref{critical euq}, we obtain
\begin{equation}
\begin{aligned}\label{1th order term}
D_p u D_{\ov p}u
\geq &\ \frac{1}{2} \frac{|D_p P_m|^2}{P_m ^2 (\phi')^2}-(\frac{\varphi^{\prime}}{\phi^{\prime}})^2 \big|D_p |Du|^2 \big|^2.
\end{aligned}
\end{equation}
Substituting \eqref{4th order term after using equ}, \eqref{derivative of Du}, \eqref{derivative of u}, \eqref{1th order term} into \eqref{diffrrentiate G second time},
\begin{equation}
\label{main inequality 1}
\begin{aligned}
0 \geq &\ -\frac{C\sum_{j}\lambda_j^{m-1}}{P_m} (1+|DDu|^2+|D\ov{D}u|^2+(1+\lambda_1)\mathcal{F}+\lambda_1) \\
&\ +\frac{\sum_j\lambda_j^{m-1}}{P_m} \Big(-\sigma_k^{p\overline{q},r\overline{s}} D_{\ov{j}} g_{r \overline{s}}{D_{j}g_{p \overline{q}}}+\sum_{\ell}{\psi_{v_l}}g_{j \overline{j}l}
+\sum_{\ell}{\psi_{\overline{v}_{\ell}}}g_{j \overline{j}\overline{\ell}}\\
&\ -\sigma_k^{p\overline{p}}(T_{pj}^au_{a \ov p \ov j}+\ov{T_{pj}^a}u_{p \ov a j}) \Big)
+ \frac{m-1}{P_m}  \sum_{j}\lambda_j^{m-2}\sigma_k^{p\overline{p}}|{D_p}g_{j \overline{j}}|^2\\
&\ + \frac{1}{P_m} \sigma_k^{p\overline{p}}\sum_{i\ne j}{\frac{\lambda_i^{m-1}-\lambda_j^{m-1}}{\lambda_i-\lambda_j}}|{D_p}{g_{i \overline{j}}}|^2 - \Big(1-\frac{\phi''}{2({\phi'})^2}\Big) \frac{|DP_m|_{\sigma}^2}{mP_m^2} \\
&\  +\frac{\varphi'}{m}\bigg(|DDu|_{\sigma\omega}^2+\frac{1}{2}|D \ov{D} u|_{\sigma\omega}^2 \bigg)
+\bigg(-\frac{\phi'}{m}\varepsilon-C\frac{\varphi'}{m} \bigg)\mathcal{F}\\
&\ +2\frac{\varphi'}{m}\operatorname{Re}\Big( \sum_{p,m}(D_{p}D_{m}uD_{\overline{p}}u+D_{p}uD_{\overline{p}}D_mu)\psi_{v_m} \Big)\\
&\ +k\frac{\phi'}{m}\psi-C\frac{\varphi'}{m}+\frac{\sigma_k^{p\overline{p}}}{m}\Big(\varphi''-\phi'' \frac{\varphi'^2}{\phi'^2} \Big) \bigg|D_p|Du|^2\bigg|^2.
\end{aligned}
\end{equation}
Here we used $\varphi' > 0$ and $\phi' < 0$ in \eqref{testfun}.

From the critical equation \eqref{critical euq}, we have
\begin{equation}
0=\frac{m}{P_m}\sum_j \lambda_j ^{m-1}D_{\ell} g_{j \ov j}+\varphi' D_{\ell}(u_p u_{\ov q} \omega^{p \ov q})+\phi' D_{\ell}u.
\end{equation}
It follows that
\begin{equation}
\begin{aligned}
0=&\; \frac{1}{P_m}\sum_j \lambda_j ^{m-1} \sum_{\ell} \psi_{v_{\ell}}D_{\ell} g_{j \ov j} \\
&\; +\frac{\varphi'}{m} \sum_{\ell} \psi_{v_{\ell}} \sum_{p}(u_{p \ell}u_{\ov p}+u_p u_{\ov p \ell})+\frac{\phi'}{m} \sum_{\ell} u_{\ell} \psi_{v_{\ell}}.
\end{aligned}
\end{equation}
Similarly, we have
\begin{equation}
\begin{aligned}
0=&\; \frac{1}{P_m}\sum_j \lambda_j ^{m-1} \sum_{\ell} \psi_{\ov{v}_{\ell}}D_{\ov{\ell}} g_{j \ov j} \\
&\; +\frac{\varphi'}{m} \sum_{\ell} \psi_{\ov{v}_{\ell}} \sum_{p}(u_{p \ov{\ell}}u_{\ov p}+u_p u_{\ov p \ov{\ell}})+\frac{\phi'}{m} \sum_{\ell} u_{\ov{\ell}} \psi_{\ov{v}_{\ell}}.
\end{aligned}
\end{equation}
Then we have
\begin{equation}
\begin{aligned}
&\frac{1}{P_m}\sum_j \lambda_j ^{m-1} \sum_{\ell} (\psi_{v_{\ell}}D_{\ell} g_{j \ov j}+\psi_{\ov{v}_{\ell}}D_{\ov{\ell}} g_{j \ov j})\\
&+2\frac{\varphi'}{m} \operatorname{Re}\bigg\{\sum_{p,m}(D_p D_m u D_{\ov p}u+D_p u D_{\ov p}D_m u)\psi_{v_m}\bigg\} \\
=&\ -\frac{\phi'}{m} \sum_{\ell}(u_{\ell} \psi_{v_{\ell}}+
u_{\ov{\ell}} \psi_{\ov{v}_{\ell}})-\frac{\varphi'}{m}\operatorname{Re}(\sum_{p,k,m}\psi_{v_m}T_{pm}^{k}u_k u_{\ov p}) \\
\geq &\ C\frac{\phi'}{m}-C\frac{\varphi'}{m}.
\end{aligned}
\end{equation}
Using \eqref{symmetric func 2th deriv}, one can obtain the well-known identity
\begin{equation}
-\sigma_{k}^{p \overline{q}, r \overline{s}} D_{j} g_{p \overline{q}} D_{\overline{j}} g_{r \overline{s}}=-\sigma_{k}^{p \overline{p}, q \overline{q}} D_{j} g_{p \overline{p}} D_{\overline{j}} g_{q \overline{q}}+\sigma_{k}^{p \overline{p}, q \overline{q}}\left|D_{j} g_{q \overline{p}}\right|^{2},
\end{equation}
where $\sigma_{k}^{p \overline{p}, q \overline{q}}=\frac{\partial}{\partial \lambda_{p}} \frac{\partial}{\partial \lambda_{q}} \sigma_{k}(\lambda)$.

We may assume that $\lambda_{1} \gg 1$. By the assumption $\varphi''-\phi''(\frac{\varphi'}{\phi'})^{2} \geq 0$ in \eqref{testfun}, the main inequality \eqref{main inequality 1} becomes
\begin{equation}
\begin{aligned}\label{main inequality 2}
0 \geq &\ \frac{1}{P_m}\sum_j\lambda_j^{m-1}\bigg(-\sigma_k^{p\overline{p},q\overline{q}} {D_{j}g_{p \overline{p}}} D_{\ov{j}} g_{q \overline{q}}+\sigma_k^{p\overline{p},q\overline{q}}|D_j g_{q \ov p}|^2 \bigg) \\
&\ -\frac{2\sigma_k^{p\overline{p}}}{P_m}\sum_j\lambda_j^{m-1}\operatorname {Re}(\ov{T_{pj}^a}u_{p \ov a j})
+\frac{\sigma_k^{p\overline{p}}}{P_m}(m-1)\sum_{j}\lambda_j^{m-2}|{D_p}g_{j \overline{j}}|^2\\
&\ +\frac{\sigma_k^{p\overline{p}}}{P_m}\sum_{i\ne j}{\frac{\lambda_i^{m-1}-\lambda_j^{m-1}}{\lambda_i-\lambda_j}}|{D_p}{g_{i \overline{j}}}|^2 -\bigg(1-\frac{\phi''}{2({\phi'})^2}\bigg)\frac{|DP_m|_{\sigma}^2}{mP_m^2}\\
&\ +\frac{\varphi'}{m}\bigg(|DDu|_{\sigma\omega}^2+\frac{1}{2}|D \ov{D} u|_{\sigma\omega}^2 \bigg) +\bigg(-\frac{\phi'}{m}\varepsilon-C\frac{\varphi'}{m}-C\bigg)\mathcal{F}\\
&\ +C\frac{\phi'}{m}-C\frac{\varphi'}{m}-C
-\frac{C}{\lambda_1}(1+|DDu|^2+|D\ov{D}u|^2).
\end{aligned}
\end{equation}
Let
\begin{equation}
\begin{aligned}
{\tilde A_j} = &\; \frac{1}{P_m}\lambda_j^{m-1} \sum_{p,q} \sigma_{k}^{p \overline{p}, q \overline{q}} D_{j} g_{p \overline{p}} D_{\overline{j}} g_{q \overline{q}},\;\;\;\;
{\tilde B_q}= \frac{1}{P_{m}} \sum_{j,p} \lambda_{j}^{m-1}\sigma_{k}^{p \overline{p}, q \overline{q}}\left|D_{j} g_{q \overline{p}}\right|^{2},
\nonumber\\
C_p=&\; \frac{m-1}{P_m} \sigma_{k}^{p \overline{p}} \sum_{j} \lambda_{j}^{m-2} \left|D_{p} g_{j \overline{j}}\right|^{2},\;\;\;\;\;\;\;
{\tilde D_p}=\frac{\sigma_{k}^{p \overline{p}}}{P_m} \sum_{j \neq i} \frac{\lambda_{i}^{m-1}-\lambda_{j}^{m-1}}{\lambda_{i}-\lambda_{j}}\left|D_{p} g_{i \overline{j}}\right|^{2},\nonumber\\
E_i=&\;\frac{m}{P^2_m}\sigma_k^{i \ov i}|\sum_p \lambda_p^{m-1}D_i g_{p \ov p}|^2,\;\;\;\;\;\;\;\;\;\;\;\;\;\;
H_p= \frac{2 \sigma_{k}^{p \overline{p}}}{P_{m}} \sum_{j,a} \lambda_{j}^{m-1}\operatorname{Re}(\ov{T_{pj}^a}u_{p \ov a j})\nonumber.
\end{aligned}
\end{equation}
Then \eqref{main inequality 2} becomes
\begin{equation}
\begin{aligned}\label{main inequalit 3}
0\geq &\ -\sum_j {\tilde A_j}+\sum_q {\tilde B_q}+\sum_p C_p+\sum_p {\tilde D_p}-\sum_p H_p \\
&\ -\Big(1-\frac{\phi^{\prime \prime}}{2\left(\phi^{\prime}\right)^{2}}\Big)\sum_i E_i
+\frac{\varphi^{\prime}}{m}\left(|D D u|_{\sigma \omega}^{2}+\frac{1}{2}|D \ov{D} u|_{\sigma \omega}^{2}\right) \\
&\ +\left(-\frac{\phi^{\prime}}{m} \varepsilon-C \frac{\varphi^{\prime}}{m}-C\right) \mathcal{F}
-C\Big(- \frac{\phi^{\prime}}{m}+\frac{\varphi^{\prime}}{m}+1\Big)\\
&\ -\frac{C}{\lambda_{1}}\left(1+|D D u|^{2}+|D \ov{D} u|^2\right).
\end{aligned}
\end{equation}

We first deal with the torsion term $H_p$. By \eqref{order}, for any $0<\beta<1$, we have
\begin{equation}
\begin{aligned}
H_p
\leq &\; \frac{2\sigma_k^{p\ov p}}{P_m} \sum_{j,a} \lambda_j^{m-1} |\ov{T_{pj}^a} D_p g_{j\ov a }| + C \sigma_k^{p\ov p}\\
\leq &\; \frac{\sigma_k^{p\ov p}}{P_m} \sum_{j,a}  \Big(\beta \lambda_j^{m-2} |D_p g_{j\ov a }|^2
 + \frac{1}{\beta} \lambda_j^{m} |T_{pj}^a|^2\Big) + C \sigma_k^{p\ov p}\\
\leq &\; \frac{\sigma_k^{p\ov p}}{P_m} \beta \sum_{a\neq j}\lambda_j^{m-2} |D_p g_{j\ov a }|^2
 +\frac{\sigma_k^{p\ov p}}{P_m} \beta \sum_{j}\lambda_j^{m-2} |D_p g_{j\ov j }|^2
 + \frac{C}{\beta} \sigma_k^{p\ov p}.
\end{aligned}
\end{equation}
By direct computation, we have
\begin{align}
 \tilde D_p = \frac{\sigma_k^{p\ov p}}{P_m} \sum_{j\neq i} \sum_{s=0}^{m-2}\lambda_i^{m-2-s}\lambda_j^s |D_p g_{i\ov j}|^2
 \geq \frac{\sigma_k^{p\ov p}}{P_m} \sum_{j\neq i}  \lambda_i^{m-2}|D_p g_{i \ov j}|^2.
\end{align}
Now we see that
\begin{equation}
\label{HCD}
\begin{aligned}
&\; - \sum_p H_p + \sum_p C_p + \sum_p \tilde D_p \\
\geq &\; (1-\beta) \sum_p \tilde D_p + (1-\frac{\beta}{m-1})\sum_p C_p - \frac{C}{\beta}\mathcal{F}.
\end{aligned}
\end{equation}
Substituting \eqref{HCD} to \eqref{main inequalit 3} yields
\begin{equation}
\begin{aligned}
\label{main inequality 4}
0\geq &\; -\sum_j \tilde{A_j}+\sum_q {\tilde B_q}
+(1 - \beta) \bigg(\sum_p C_p + \sum_p {\tilde D_p} \bigg)\\
&\;-\left(1-\frac{\phi^{\prime \prime}}{2\left(\phi^{\prime}\right)^{2}}\right)\sum_i E_i+\frac{\varphi^{\prime}}{m}\left(|D D u|_{\sigma \omega}^{2}+\frac{1}{2}|D \ov{D} u|_{\sigma \omega}^{2}\right) \\
&\ +\left(-\frac{\phi^{\prime}}{m} \varepsilon-C \frac{\varphi^{\prime}}{m}-\frac{C}{\beta}\right) \mathcal{F}-C(\varphi'- \phi'+1)\\
&\ -\frac{C}{\lambda_{1}}\left(1+|D D u|^{2}+|D \ov{D} u|^2\right).
\end{aligned}
\end{equation}

We need a lemma from \cite{G-R-W}.
\begin{lemma}[\cite{G-R-W}]
Suppose $1 \leq \ell<k \leq n$, and let $\alpha=1 /(k-\ell)$. Let $W=(w_{p\overline{q}})$ be a Hermitian tensor in the $\Gamma_{k}$ cone. Then for any $\theta>0$,
\begin{equation}
\begin{aligned}
&\ -\sigma_{k}^{p \ov{p}, q \ov{q}}(W) w_{\overline{p} p i} w_{\overline{q} q \overline{i}}+\left(1-\alpha+\frac{\alpha}{\theta}\right) \frac{\left|D_{i} \sigma_{k}(W)\right|^{2}}{\sigma_{k}(W)}\\
\geq &\ \sigma_{k}(W)(\alpha+1-\alpha \theta)\left|\frac{D_{i} \sigma_{\ell}(W)}{\sigma_{\ell}(W)}\right|^{2}-\frac{\sigma_{k}}{\sigma_{\ell}}(W) \sigma_{\ell}^{p \overline{p}, q \overline{q}}(W) w_{\overline{p} p i} w_{\overline{q} q \overline{i}}.
\end{aligned}
\end{equation}
\end{lemma}
Taking $\ell=1$ in the above lemma, we have
\begin{align}
-\sigma_k^{p\ov p, q\ov q} D_i g_{p\ov p}D_{\ov i} g_{q\ov q} + K |D_i \sigma_k|^2 \geq 0,
\end{align}
for $K>(1-\alpha+\alpha / \theta)(\inf \psi)^{-1}$ if $2 \leq k\leq n$.
We shall denote
\[\begin{aligned}
A_j = &\; \frac{1}{P_m} \lambda_j^{m-1} \Big(K |D_j \sigma_k|^2 - \sigma_k^{p\ov p, q\ov q} D_j g_{p\ov p}D_{\ov j} g_{q\ov q}\Big),\\
B_q = &\; \frac{1}{P_m} \sum_p \lambda_p^{m-1} \sigma_k^{p\ov p, q\ov q}|D_q g_{p\ov p}|^2,\\
D_i = &\; \frac{1}{P_m} \sum_{p\neq i} \sigma_k^{p\ov p} \frac{\lambda_p^{m-1} - \lambda_i^{m-1}}{\lambda_p -\lambda_i} |D_i g_{p\ov p}|^2.
\end{aligned}\]
Define $H_{j \ov p q} = D_j\chi_{q\ov p} - D_q \chi_{j\ov p}$. For any constant $0<\tau <1$, we can estimate
\begin{equation}\nonumber
\begin{aligned}
\sum_q \tilde B_q
\geq &\;
\frac{1}{P_m} \sum_{j,q} \lambda_j^{m-1} \sigma_k^{j\ov j, q\ov q} |D_j g_{q \ov j }|^2\\
= &\; \frac{1}{P_m} \sum_{j,q} \lambda_j^{m-1} \sigma_k^{j\ov j, q\ov q} |D_q u_{j\ov j} - T_{jq}^a u_{a \ov j} + D_q \chi_{j\ov j} + H_{j\ov j q}|^2\\
\geq &\; \frac{1}{P_m}\sum_{q,j} \lambda_j^{m-1} \sigma_k^{j\ov j, q\ov q}
\Big( (1-\tau)|D_q g_{j \ov j}|^2 - (\frac{1}{\tau} - 1)| H_{j\ov j q}- T_{jq}^a u_{a\ov j} |^2\Big)\\
= &\; (1-\tau) \sum_q B_q - \frac{\frac{1}{\tau}-1}{P_m} \sum_{q,j} \lambda_j^{m-2}\big(\sigma_k^{j\ov j, q\ov q} \lambda_j\big)
| H_{j\ov j q} - T_{jq}^a u_{a\ov j} |^2.
\end{aligned}
\end{equation}
Now we use $\sigma_l (\lambda|i)$ and $\sigma_l(\lambda|ij)$ to denote the $l$-th elementary function of
\[
(\lambda|i) = (\lambda_1, \ldots, \hat \lambda_i, \ldots, \lambda_n) \in \mathbb{R}^{n-1}
\]
and
\[(\lambda|ij) = (\lambda_1, \ldots, \hat\lambda_i, \ldots, \hat\lambda_j, \ldots, \lambda_n)\in \mathbb{R}^{n-2}
\]
Then, we have $\sigma_k^{i\ov i} =\sigma_{k-1} (\lambda|i)$, $\sigma_k^{p\ov p, i\ov i} = \sigma_{k-2} (\lambda|pi)$.
By the formula $\sigma_l(\lambda) = \sigma_l(\lambda|p) + \lambda_p \sigma_{l-1}(\lambda| p)$ for any $1\leq p\leq n$, we obtain
\[
\frac{1}{\tau P_m} \sum_{q,j} \lambda_j^{m-2}\big(\sigma_k^{q\ov q} - \sigma_{k-1}(\lambda|jq)\big) |H_{j\ov j q}- T_{jq}^a u_{a\ov j} |^2
\leq \frac{C}{\tau}\mathcal{F},
\]
which implies
\begin{align}\nonumber
\sum_q \tilde B_q
\geq &\; (1-\tau)\sum_q B_q -\frac{C}{\tau}\mathcal{F} \nonumber.
\end{align}
Similarly, we may estimate
\[\begin{aligned}
\sum_p \tilde D_p \geq &\; \frac{1}{P_m} \sum_{j\neq i} \sigma_k^{j\ov j}
\frac{\lambda_i^{m-1}- \lambda_j^{m-1}}{\lambda_i -\lambda_j}|D_j g_{i\ov j}|^2\\
\geq &\; \frac{1}{P_m} \sum_{j\neq i} \sigma_k^{j\ov j} \frac{\lambda_i^{m-1}- \lambda_j^{m-1}}{\lambda_i -\lambda_j}
\Big((1-\tau) |D_i g_{j\ov j}|^2 - \frac{C}{\tau} \lambda_1^2 \Big)\\
\geq &\; (1-\tau)\sum_i D_i - \frac{C}{\tau} \mathcal{F}.
\end{aligned}\]
Note that $\frac{\lambda_j^{m-1}}{P_m} |D_j\sigma_k|^2 \leq \frac{C}{\lambda_1} (|DDu|^2+ |D\ov D u|^2)$.
Then \eqref{main inequality 4} becomes
\begin{equation}
\begin{aligned}
\label{main inequality 5}
0\geq &\; \sum_i A_i+ (1-\tau)\sum_i B_i
+ (1-\beta) \sum_i C_i \\
&\;+(1-\beta)(1-\tau)\sum_i D_i -\left(1-\frac{\phi^{\prime \prime}}{2\left(\phi^{\prime}\right)^{2}}\right)\sum_i E_i \\
&\ +\frac{\varphi^{\prime}}{m}\left(|D D u|_{\sigma \omega}^{2}+\frac{1}{2}|D \ov{D} u|_{\sigma \omega}^{2}\right)
-\frac{C(K)}{\lambda_{1}}\left(|D D u|^{2}+|D \ov{D} u|^2\right)\\
&\ +\left(-\frac{\phi^{\prime}}{m} \varepsilon-C \frac{\varphi^{\prime}}{m}-\frac{C}{\beta} - \frac{C}{\tau}\right) \mathcal{F}
-C(\varphi'- \phi'+1),
\end{aligned}
\end{equation}
when $\lambda_1$ is sufficiently large.
Let $1-\delta = (1-\beta)(1-\tau)$. We then have
\begin{equation}
\begin{aligned}
\label{main inequality 6}
0\geq &\; (1-\delta)\sum_i \Big(  A_i+  B_i+ C_i+  D_i\Big)-\left(1-\frac{\phi^{\prime \prime}}{2\left(\phi^{\prime}\right)^{2}}\right)\sum_i E_i\\
&\; +\frac{\varphi^{\prime}}{m}\left(|D D u|_{\sigma \omega}^{2}+\frac{1}{2}|D \ov{D} u|_{\sigma \omega}^{2}\right)
 -\frac{C(K)}{\lambda_{1}}\left(|D D u|^{2}+|D \ov{D} u|^2\right)\\
&\ +\left(-\frac{\phi^{\prime}}{m} \varepsilon-C \frac{\varphi^{\prime}}{m}-\frac{C}{\beta} - \frac{C}{\tau}\right) \mathcal{F}
-C(\varphi'- \phi'+1).
\end{aligned}
\end{equation}

Now we choose $\phi$ and $\varphi$ to satisfy \eqref{testfun}. Let $\varphi (t) = e^{Nt}$ and $\phi(s) = e^{M(-s+L)}$ where $L \geq |u|_{C^1}+1$ is a constant.
Then, we see
\[
\varphi''-\phi'' \frac{\varphi'^2}{\phi'^2} = N^2 e^{Nt} - \frac{N^2 e^{2Nt}}{e^{M(-s +L)}} >0, \quad \varphi' >0, \quad \phi' < 0,
\]
when $M \gg N > 1$, which shows the assumption \eqref{testfun} is satisfied. Since, at $p$,
\[
\frac{\phi''}{2(\phi')^2} = \frac{1}{2e^{M(-u(p)+L)}},
\]
by choosing
\[
\beta =\tau = \frac{1}{6 e^{M(-u(p) + L)}},
\]
we obtain that $1-\delta \geq 1-\frac{\phi^{\prime \prime}}{2\left(\phi^{\prime}\right)^{2}} $.
By the arguments as in \cite{P-P-Z C^2}, we may assume without loss of generality that
\[
A_{i}+B_{i}+C_{i}+D_{i}-E_{i} \geq 0, \quad \forall i=1, \ldots, n.
\]
Since $\sigma_{k}^{i \overline{i}} \geq \sigma_{k}^{1 \overline{1}} \geq \frac{k}{n} \frac{\sigma_{k}}{\lambda_{1}} \geq \frac{1}{C \lambda_{1}}$ for fixed $i$, we can estimate
\[
|D D u|_{\sigma \omega}^{2}+\frac{1}{2}|D \ov{D} u|_{\sigma \omega}^{2} \geq \frac{1}{C \lambda_{1}}\left(|D D u|^{2}+|D \overline{D} u|^{2}\right) \geq \frac{1}{C\lambda_1} |DDu|^2 + \frac{\lambda_1}{C}.
\]
Now \eqref{main inequality 6} becomes
\begin{equation}
\begin{aligned}
0\geq &\; \Big(\frac{\varphi'}{m C} - C(K)\Big)\lambda_1 + \frac{1}{\lambda_1} \Big(\frac{\varphi'}{m C} - C(K)\Big) |DDu|^2 \\
&\; +\left(-\frac{\phi^{\prime}}{m} \varepsilon-C \frac{\varphi^{\prime}}{m}-\frac{C}{\beta} - \frac{C}{\tau}\right) \mathcal{F} -C(\varphi'- \phi'+1).
\end{aligned}
\end{equation}
Taking $N$ large enough, we can ensure that $\frac{\varphi'}{m C} - C(K)>0$.
For fixed $N$, it follows that
\[-\frac{\phi^{\prime}}{m} \varepsilon-C \frac{\varphi^{\prime}}{m}-\frac{C}{\beta} - \frac{C}{\tau}
=\frac{M}{m} \varepsilon \phi - C \frac{N}{m} \varphi - C \phi > 0 \]
when $M\gg N$.
This leads to
\[\begin{aligned}
0\geq &\; \Big(\frac{\varphi'}{m C} - C(K)\Big)\lambda_1  - C,
\end{aligned}\]
which finally gives an upper bound of $\lambda_1$.

\section{ Outline of proof of Theorem \ref{thm 2} }
We can rewrite (\ref{Sk equation 2}) as follows:
\begin{equation}\label{Rewrite Sk equ 2}
{\sigma}_{k}(\tilde g)={\sigma}_{k}\left(\chi_{i \overline{j}}+u_{i \overline{j}}+a_{i} u_{\overline{j}}+a_{\overline{j}} u_{i}\right)=\psi(z, D u, u).
\end{equation}
By direct calculation and Cauchy inequality, we have
\begin{equation}
\begin{aligned}
&\ { \sigma}_{k}^{p \overline{q}} D_{\overline{q}} D_{p} \tilde g_{i \overline{j}}\\
\geq &\ -{ \sigma}_{k}^{p \overline{q}, r \overline{s}} D_{\ov{j}} \tilde g_{r \overline{s}} D_{i} \tilde g_{p \overline{q}}+\sum_\ell \psi_{v_\ell} \tilde g_{i \ov{j} \ell} + \sum_\ell \psi_{\ov{v}_\ell} \tilde g_{i \ov{j} \ov{\ell}}-{ \sigma}_{k}^{p \ov{q}} (T_{p i}^{s} u_{s \overline{q} \overline{j}} \\
&\ +\overline{T_{q j}^{s}} u_{p \overline{s} i})-C(1+|DDu|^{2}+|D \ov{D}u|^{2}+\mathcal{F}+\lambda_{1} \mathcal{ F}+\lambda_{1})\\
&\ +{ \sigma}_{k}^{p \overline{q}}(a_i u_{\ov j p \ov q}+a_{\ov j} u_{i p \ov q}-a_p u_{\ov q i \ov j}-a_{\ov q} u_{pi \ov j})
- |D \ov{D} u|_{{ \sigma} \omega}^{2}- |DDu|_{{\sigma} \omega}^{2}.
\end{aligned}
\end{equation}
and
\begin{equation}
\begin{aligned}
{ \sigma}_{k}^{p \overline{q}}|D u|_{p \overline{q}}^{2}
\geq &\; 2 \operatorname{Re}\left\{\sum_{p, m}\left(D_{p} D_{m} u D_{\overline{p}} u+D_{p} u D_{\overline{p}} D_{m} u\right) \psi_{v_{m}}\right\}\\
&\; -C-C\mathcal{F}+\frac{1}{2}|D D u|_{{ \sigma} \omega}^{2}+\frac{1}{2}|D \ov{D} u|_{{ \sigma} \omega}^{2}.
\end{aligned}
\end{equation}
We also have
\begin{equation}
\begin{aligned}
&\ \Big|\frac{2{\sigma}_{k}^{p \overline{p}}}{m} \operatorname{Re} \left\{ a_p (\frac{D_{\ov p} {{P}_m}}{{P}_m}+\varphi' D_{\ov p}|Du|^2) \right\} \Big| \\
\leq &\ \frac{\phi''}{4{\phi'}^2} \frac{|D{P_m}|_{{\sigma}}^2}{m{P}^2_m}+ \frac{{ \sigma}_{k}^{p \overline{p}}}{m} {\varphi'}^2 \frac{\phi''}{\phi'^2}\bigg|D_p|Du|^2\bigg|^2+C\frac{\phi'^2}{\phi''} \mathcal{ F}.
\end{aligned}
\end{equation}

Now we apply the maximum principle to the test function \eqref{test} bearing in mind that
$\lambda_{1}, \lambda_{2}, \dots, \lambda_{n}$ are the eigenvalues of $\tilde g_{i \overline{j}} =\chi_{i \overline{j}}+u_{i \overline{j}}+a_{i} u_{\overline{j}}+a_{\overline{j}} u_{i}$ with respect to $\omega$.
Instead of the assumptions in \eqref{testfun}, we now require
\begin{equation}
\label{testfun'}
\varphi''-2\phi''(\frac{\varphi'}{\phi'})^{2} \geq 0, \;\varphi' > 0, \; \phi' < 0.
\end{equation}
With the above calculations, similar to \eqref{main inequality 2}, we have
\begin{equation}
\begin{aligned}\label{main inequality 2'}
0 \geq
&\ \frac{1}{{P_m}}\sum_j{{\lambda}}_j^{m-1}\bigg(-{\sigma}_k^{p\overline{p},q\overline{q}} {D_{j}\tilde g_{p \overline{p}}} D_{\ov{j}} \tilde g_{q \overline{q}}+{ \sigma}_k^{p\overline{p},q\overline{q}}|D_j \tilde g_{q \ov p}|^2 \bigg) \\
&\ -\frac{2{ \sigma}_k^{p\overline{p}}}{P_m}\sum_j{\lambda}_j^{m-1}
\operatorname{Re}\Big(\ov{T_{pj}^s}u_{p \ov s j}\Big)
+ \frac{(m-1)\sigma_k^{p\overline{p}}}{P_m} \sum_j \lambda_j^{m-2}|{D_p} \tilde g_{j \overline{j}}|^2 \\
&\ +\frac{{ \sigma}_k^{p\overline{p}}}{{{P}_m}} \sum_{i\ne j}{\frac{\lambda_i^{m-1}
- \lambda_j^{m-1}}{\lambda_i-\lambda_j}}|{D_p}{\tilde g_{i \overline{j}}}|^2
-\bigg(1-\frac{\phi''}{4({\phi'})^2}\bigg)\frac{ |DP_m|_{\sigma}^2 }{m P_m^2} \\
&\ +\Big(\frac{\varphi'}{2 m} - \frac{1}{\lambda_1}\Big) \bigg(|DDu|_{{ \sigma}\omega}^2+|D \ov{D} u|_{{ \sigma}\omega}^2 \bigg)
+C\frac{\phi'}{m}-C\frac{\varphi'}{m}-C\\
&\ +\bigg(-\frac{\phi'}{m}\varepsilon-C\frac{\varphi'}{m}-C-C\frac{\phi'^2}{\phi''}\bigg)\mathcal{ F}
-\frac{C}{\lambda_1}(|DDu|^2+|D\ov{D}u|^2)\\
&\  + \frac{{ \sigma}_k^{p\overline{p}}}{{P_m}}\sum_j \lambda_j^{m-1}
\Big(a_j u_{\ov j p \ov p}+a_{\ov j} u_{j p \ov p}-a_p u_{\ov p j \ov j}-a_{\ov p} u_{pj \ov j} \Big) ,
\end{aligned}
\end{equation}
where we used $\varphi''- 2 \phi'' \frac{\varphi'^2}{\phi'^2} \geq 0$ in \eqref{testfun'} and $\lambda_{1} \gg 1$.

We use the same notation for third order terms as in Section 3.
Denote
\[
I_p= \frac{{\sigma}_k^{p\overline{p}}}{P_m}\sum_j \lambda_j^{m-1}\Big(a_j u_{\ov j p \ov p}+a_{\ov j} u_{j p \ov p}-a_p u_{\ov p j \ov j}-a_{\ov p} u_{pj \ov j}\Big).
\]
For any $0<\beta<1$, we can estimate
\begin{equation}\nonumber
\begin{aligned}
I_p
\leq &\; \frac{\beta}{4} \frac{{\sigma}_k^{p\ov p}}{{{P}_m}} \sum_{j}  {{\lambda}}_j^{m-2}
\Big(|D_p \tilde g_{j\ov p }|^2 + |D_p \tilde g_{j\ov j }|^2\Big)
+ \frac{C}{\beta} {\sigma}_k^{p\ov p}+\frac{C}{{{\lambda}}_1} { \sigma}_k^{p\ov p} \sum_{j} |u_{jp}|^2\\
\leq &\; \frac{\beta}{2} \frac{{ \sigma}_k^{p\ov p}}{{{P}_m}} \sum_{j,s}  {{\lambda}}_j^{m-2} |D_p \tilde g_{j\ov s }|^2+ \frac{C}{\beta} {\sigma}_k^{p\ov p}+\frac{C}{{{\lambda}}_1} {\sigma}_k^{p\ov p} \sum_{j} |u_{jp}|^2
\end{aligned}
\end{equation}
and
\begin{equation}\nonumber
\begin{aligned}
H_p
\leq &\; \frac{2{ \sigma}_k^{p\ov p}}{{{P}_m}} \sum_{j,s} {{\lambda}}_j^{m-1} |\ov{T_{pj}^s} D_p \tilde g_{j\ov s }| + C { \sigma}_k^{p\ov p}+\frac{C}{{{\lambda}}_1} { \sigma}_k^{p\ov p}\sum_{j} |u_{jp}|^2\\
\leq &\;  \frac{\beta}{2} \frac{{ \sigma}_k^{p\ov p}}{{P}_m} \sum_{j,s} {{\lambda}}_j^{m-2} |D_p \tilde g_{j\ov s }|^2+ \frac{C}{\beta} { \sigma}_k^{p\ov p}+\frac{C}{{{\lambda}}_1} {\sigma}_k^{p\ov p}\sum_{j} |u_{jp}|^2.
\end{aligned}
\end{equation}
Now we have
\begin{equation}\nonumber
\begin{aligned}
&\; - \sum_p H_p -\sum_p I_p+ \sum_p C_p + \sum_p \tilde D_p \\
\geq &\; (1-\beta) \sum_p \tilde D_p + (1-\beta)\sum_p C_p - \frac{C}{\beta}\mathcal{ F}-\frac{C}{\lambda_1}|DDu|_{{ \sigma} \omega}^2,
\end{aligned}
\end{equation}
Define $H'_{j \ov q p} = D_j(\chi_{p\ov q}+a_p u_{\ov q}+a_{\ov q} u_p) - D_p (\chi_{j\ov q}+a_j u_{\ov q}+a_{\ov q} u_j)$.
Note that $|H_{j\ov j q}'|\leq C+C\lambda_1$.
For any $0<\tau <1$, we can estimate
\begin{align}\nonumber
\sum_q \tilde B_q +\sum_p \tilde D_p
\geq &\; (1-\tau)\Big(\sum_q B_q + \sum_i D_i\Big)-\frac{C}{\tau}\mathcal{F}.
\end{align}
Let $1- \delta = (1-\beta)(1-\tau)$. Similar to \eqref{main inequality 6},
when $\frac{C}{\lambda_1} \leq \frac{\varphi'}{4m}$, we obtain
\begin{equation}
\begin{aligned}
\label{main inequality 6'}
0\geq &\; (1-\delta)\sum_i \Big(  A_i+  B_i+ C_i+  D_i\Big)
-\left(1-\frac{\phi^{\prime \prime}}{4\left(\phi^{\prime}\right)^{2}}\right)\sum_i E_i\\
&\; + \frac{\varphi'}{4 m} \left(|D D u|_{{ \sigma} \omega}^{2}+|D \ov{D} u|_{{ \sigma} \omega}^{2}\right)
 -\frac{C(K)}{\lambda_{1}}\left(|D D u|^{2}+|D \ov{D} u|^2\right)\\
&\ +\left(-\frac{\phi^{\prime}}{m} \varepsilon-C \frac{\varphi^{\prime}}{m}-\frac{C}{\beta}
- \frac{C}{\tau}-C\frac{\phi'^2}{\phi''}\right) \mathcal{ F}
-C(\varphi'- \phi'+1).
\end{aligned}
\end{equation}

Still choose $\phi$ and $\varphi$ as before and \eqref{testfun'} is guaranteed by $M\gg N$.
If we choose $\beta=\tau = \frac{\phi''}{12(\phi')^2}$, we get $1-\delta \geq 1-\frac{\phi^{\prime \prime}}{ 4 \left(\phi^{\prime}\right)^{2}}$.
Now the proof is the same as that of Theorem \ref{thm 1}.

\end{document}